\PassOptionsToPackage{unicode}{hyperref}
\PassOptionsToPackage{hyphens}{url}
\documentclass[a4paper]
{article}
\usepackage{authblk}
\usepackage{amsmath,amssymb}
\usepackage{lmodern}
\usepackage{iftex}
\usepackage[top=2cm, bottom=2cm, left=2cm, right=2cm]{geometry}

\ifPDFTeX
  \usepackage[T1]{fontenc}
  \usepackage[utf8]{inputenc}
  \usepackage{textcomp} 
\else 
  \usepackage{unicode-math}
  \defaultfontfeatures{Scale=MatchLowercase}
  \defaultfontfeatures[\rmfamily]{Ligatures=TeX,Scale=1}
\fi
\IfFileExists{upquote.sty}{\usepackage{upquote}}{}
\IfFileExists{microtype.sty}{
  \usepackage[]{microtype}
  \UseMicrotypeSet[protrusion]{basicmath} 
}{}
\makeatletter
\@ifundefined{KOMAClassName}{
  \IfFileExists{parskip.sty}{%
    \usepackage{parskip}
  }{
    \setlength{\parindent}{0pt}
    \setlength{\parskip}{6pt plus 2pt minus 1pt}}
}{
  \KOMAoptions{parskip=half}}
\makeatother
\usepackage{xcolor}
\IfFileExists{xurl.sty}{\usepackage{xurl}}{} 
\IfFileExists{bookmark.sty}{\usepackage{bookmark}}{\usepackage{hyperref}}
\hypersetup{
  pdftitle={FlowClass.jl: Classifying Dynamical Systems by Structural Properties in Julia},
  hidelinks,
  pdfcreator={LaTeX via pandoc}}
\usepackage{color}
\usepackage{fancyvrb}

\DefineVerbatimEnvironment{Highlighting}{Verbatim}{commandchars=\\\{\}}
\newenvironment{Shaded}{}{}

\newcommand{\CommentTok}[1]{\textcolor[rgb]{0.38,0.63,0.69}{\textit{#1}}}

\newcommand{\DataTypeTok}[1]{\textcolor[rgb]{0.56,0.13,0.00}{#1}}

\newcommand{\FloatTok}[1]{\textcolor[rgb]{0.25,0.63,0.44}{#1}}

\newcommand{\KeywordTok}[1]{\textcolor[rgb]{0.00,0.44,0.13}{\textbf{#1}}}
\newcommand{\NormalTok}[1]{#1}
\newcommand{\OperatorTok}[1]{\textcolor[rgb]{0.40,0.40,0.40}{#1}}

\newcommand{\PreprocessorTok}[1]{\textcolor[rgb]{0.74,0.48,0.00}{#1}}

\newcommand{\SpecialCharTok}[1]{\textcolor[rgb]{0.25,0.44,0.63}{#1}}

\newcommand{\StringTok}[1]{\textcolor[rgb]{0.25,0.44,0.63}{#1}}

\usepackage{longtable,booktabs,array}
\usepackage{calc} 
\usepackage{etoolbox}
\makeatletter
\patchcmd\longtable{\par}{\if@noskipsec\mbox{}\fi\par}{}{}
\makeatother
\IfFileExists{footnotehyper.sty}{\usepackage{footnotehyper}}{\usepackage{footnote}}
\makesavenoteenv{longtable}
\usepackage{graphicx}
\makeatletter
\def\maxwidth{\ifdim\Gin@nat@width>\linewidth\linewidth\else\Gin@nat@width\fi}
\def\maxheight{\ifdim\Gin@nat@height>\textheight\textheight\else\Gin@nat@height\fi}
\makeatother
\setkeys{Gin}{width=\maxwidth,height=\maxheight,keepaspectratio}
\makeatletter
\def\fps@figure{htbp}
\makeatother
\providecommand{\tightlist}{%
  \setlength{\itemsep}{0pt}\setlength{\parskip}{0pt}}
\ifLuaTeX
  \usepackage{selnolig}  
\fi

\title{FlowClass.jl: Classifying Dynamical Systems by Structural
Properties in Julia}
\author[1,2,3]{Michael P.H. Stumpf}
\affil[1]{School of BioScience, University of Melbourne, Melbourne, Australia}
\affil[2]{School of Mathematics and Statistics, University of Melbourne, Melbourne, Australia}
\affil[3]{Cell Bauhaus PTY LTY, University of Melbourne, Melbourne, Australia}

\date{15 December 2025}

\begin{document}
\maketitle

\hypertarget{summary}{%
\section{Summary}\label{summary}}

\texttt{FlowClass.jl} is a Julia package
\cite{Bezanson:2017,Roesch2023} for classifying continuous-time
dynamical systems into a hierarchy of structural classes: Gradient,
Gradient-like, Morse-Smale, Structurally Stable, and General. Given a
vector field \(\mathbf{F}(\mathbf{x})\) defining the system
\(\mathrm{d}\mathbf{x}/\mathrm{d}t = \mathbf{F}(\mathbf{x})\), the
package performs a battery of computational tests---Jacobian symmetry
analysis, curl magnitude estimation, fixed point detection and stability
classification, periodic orbit detection, and stable/unstable manifold
computation---to determine where the system sits within the
classification hierarchy. This classification has direct implications
for qualitative behaviour: gradient systems cannot oscillate,
Morse-Smale systems are structurally stable in less than 3 dimensions,
and general systems may exhibit chaos. Much of classical developmental
theory going back to Waddington's epigenetic landscape
\cite{Waddington1957} rests on an implicit assumption of gradient
dynamics.

The package is designed with applications in systems and developmental
biology in mind, particularly the analysis of gene regulatory networks
and cell fate decision models in the context of Waddington's epigenetic
landscape. It provides tools to assess whether a landscape metaphor is
appropriate for a given dynamical model, and to quantify the magnitude
of non-gradient (curl) dynamics.

\hypertarget{statement-of-need}{%
\section{Statement of Need}\label{statement-of-need}}

Dynamical systems models are ubiquitous in biology, physics, and
engineering \cite{Strogatz2015}. A central question when analysing such
models is whether the system can be understood in terms of a potential
landscape---that is, whether trajectories follow the gradient of a
potential function toward attracting states. This question is especially
pertinent in developmental biology \cite{Huang2007,Moris2016}, where
Waddington's epigenetic landscape provides a powerful metaphor for cell
differentiation \cite{Waddington1957}. However, as recent work has
emphasised, most biological systems are \emph{not} gradient systems:
they exhibit curl dynamics arising from non-reciprocal interactions in
gene regulatory networks \cite{Brackston2018,Wang2015}.

Despite the theoretical importance of this distinction
\cite{Smale1967,PalisDeMelo1982}, there has been no comprehensive
software package for computationally classifying dynamical systems by
their structural properties. Existing tools focus on specific
aspects---bifurcation analysis, Lyapunov exponent computation, or
trajectory simulation---but do not provide an integrated framework for
structural classification. \texttt{FlowClass.jl} fills this gap by
implementing a systematic classification pipeline that moves from the
most restrictive class (gradient systems) to the most general, providing
quantitative measures at each stage. \includegraphics{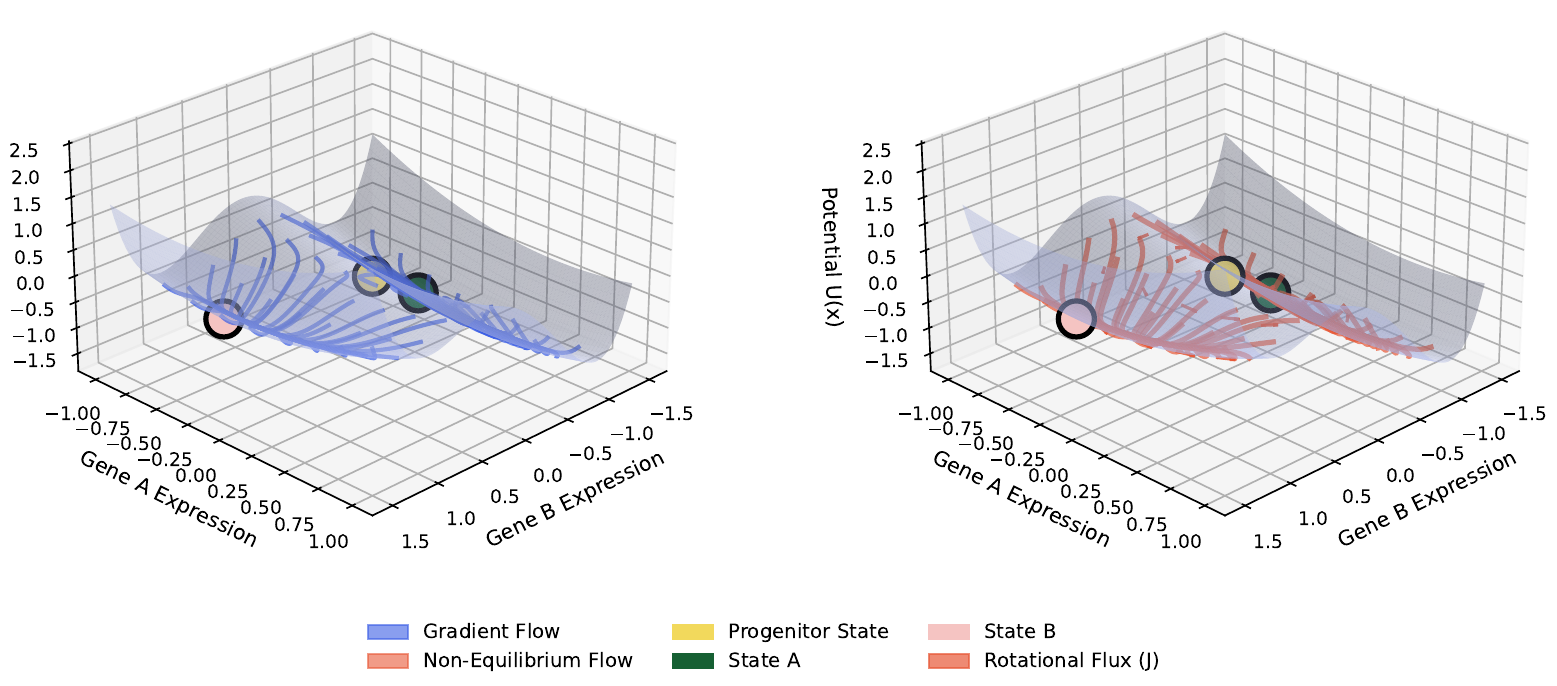}
\emph{Figure 1. Waddington's epigenetic landscape showing potential
wells corresponding to cell fates. On the left hand side the flow is
purely gradient. On the right hand side the flow has curl components and
the dynamics are no longer determined by the landscape.}

The classification has practical consequences. For gradient systems, the
potential function fully characterises the dynamics, and forward and
reverse transition paths coincide. For systems with even moderate curl
fluxes, minimum action paths for transitions in opposite directions
differ, which has implications for understanding differentiation versus
reprogramming in stem cell biology \cite{Brackston2018,Guillemin2021}.
By quantifying the curl-to-gradient ratio, researchers can assess the
validity of landscape-based analyses for their specific models. This for
example, makes it possible to infer aspects of stem cell dynamics from
single cell data \cite{Liu2024}.

\hypertarget{classification-of-dynamical-systems}{%
\section{Classification of Dynamical
Systems}\label{classification-of-dynamical-systems}}

\textbf{Table 1.} Properties of dynamical system classes in the
classification hierarchy.

\begin{longtable}[]{@{}
  >{\raggedright\arraybackslash}p{(\columnwidth - 10\tabcolsep) * \real{0.13}}
  >{\raggedright\arraybackslash}p{(\columnwidth - 10\tabcolsep) * \real{0.18}}
  >{\raggedright\arraybackslash}p{(\columnwidth - 10\tabcolsep) * \real{0.17}}
  >{\raggedright\arraybackslash}p{(\columnwidth - 10\tabcolsep) * \real{0.13}}
  >{\raggedright\arraybackslash}p{(\columnwidth - 10\tabcolsep) * \real{0.29}}
  >{\raggedright\arraybackslash}p{(\columnwidth - 10\tabcolsep) * \real{0.10}}@{}}
\toprule
Property & Gradient & Gradient-like & Morse-Smale & Generic & General \\
\midrule
\endhead
Curl & Zero everywhere & Near zero/small & Non-zero possible & Any value
& Any value \\
Jacobian & Symmetric $(J = J^T)$ & Nearly symmetric & No requirement & No
requirement & No requirement \\
Path integrals & Path-independent & Nearly path-independent &
Path-dependent & Path-dependent & Path-dependent \\
Periodic orbits & None & None & Hyperbolic only & Hyperbolic \&
non-hyperbolic & Any \\
Critical points & All hyperbolic & All hyperbolic & All hyperbolic &
Non-hyperbolic possible at bifurcations & Any \\
Lyapunov function & Global (potential V) & Global & Local only & Local
(away from bifurcations) & May not exist \\
Transversality & N/A (no periodic orbits) & N/A & Required & May have
tangencies & Not required \\
\bottomrule
\end{longtable}

Our ability to make qualitative statements about a dynamical system
depends crucially on the nature of the dynamics
\cite{Smale1967,PalisDeMelo1982,BrackstonPhysRevE2018}. The
qualitative aspects can deliver profound biological systems. If stem
cell differentiation were to follow gradient dynamics, for example, then
the forward and backward paths through gene expression space would be
identical \cite{Guillemin2021,Vittadello2025}. For gradient systems,
and gradient-like systems we have access to Lyapunov functions, and
concepts from catastrophe theory can be applied and yield powerful
insights \cite{Rand2021}.

For more general dynamical systems we cannot fall back on such elegant
theory. A main focus in the design of\texttt{FlowClass.jl} was to
classify dynamical systems into the relevant categories that determine
whether or not a given system (at least in the specific parameterisation
considered). The routines provided as part of the package make it
possible to determine the characterising features of different classes
of dynamical systems, (c.f. Table 1).

\hypertarget{key-features}{%
\section{Key Features}\label{key-features}}

\textbf{Jacobian analysis.} The package computes Jacobian matrices using
automatic differentiation via \texttt{ForwardDiff.jl} and tests for
symmetry. For a gradient system \(\mathbf{F} = -\nabla V\), the Jacobian
equals the negative Hessian and is therefore symmetric. The relative
symmetry error \(\|(J - J^\top)/2\|_F / \|J\|_F\) provides a
scale-independent measure of deviation from gradient structure.

\textbf{Curl quantification.} For 2D and 3D systems,
\texttt{FlowClass.jl} computes the curl directly. For higher-dimensional
systems, it uses the Frobenius norm of the antisymmetric part of the
Jacobian as a generalised curl measure. The curl-to-gradient ratio
indicates the relative strength of rotational versus potential-driven
dynamics.

\textbf{Fixed point analysis.} Multi-start optimisation via
\texttt{NLsolve.jl} locates fixed points within user-specified bounds.
Each fixed point is classified by eigenvalue analysis into stable nodes,
unstable nodes, saddles, foci, centres, or non-hyperbolic points. The
presence of non-hyperbolic fixed points excludes the system from the
Morse-Smale class.

\textbf{Periodic orbit detection.} Poincaré section methods detect limit
cycles, and Floquet multiplier analysis determines their stability.
Gradient and gradient-like systems cannot possess periodic orbits; their
presence indicates at least Morse-Smale classification.

\textbf{Manifold computation.} For saddle points, the package computes
stable and unstable manifolds by integrating trajectories from initial
conditions perturbed along eigenvector directions. Transversality of
manifold intersections is assessed numerically as this is a requirement
for Morse-Smale structure.

\textbf{Integrated classification.} The \texttt{classify\_system}
function orchestrates all analyses and returns a
\texttt{ClassificationResult} containing the system class, all detected
invariant sets, quantitative measures (Jacobian symmetry, curl ratio),
and a confidence score.

\hypertarget{example-stem-cell-differentiation}{%
\section{Example: Stem Cell
Differentiation}\label{example-stem-cell-differentiation}}

As a demonstration, \texttt{FlowClass.jl} includes an implementation of
the stem cell differentiation model from \cite{Brackston2018}, who built on
earlier work by \cite{Chickarmane2012}, which describes the dynamics of
pluripotency factors (Nanog, Oct4-Sox2, Fgf4) and the differentiation
marker Gata6. The model exhibits multiple stable states corresponding to
pluripotent and differentiated cell fates, with a saddle point
representing the transition state. Classification reveals significant
curl dynamics, confirming that the system is not gradient and that
differentiation and reprogramming paths will differ---a key biological
insight.

\begin{Shaded}
\begin{Highlighting}[]
\KeywordTok{using}\NormalTok{ FlowClass}

\NormalTok{ds }\OperatorTok{=}\NormalTok{ DynamicalSystem(stem\_cell\_model}\OperatorTok{,} \FloatTok{4}\NormalTok{)}
\NormalTok{bounds }\OperatorTok{=}\NormalTok{ ((}\FloatTok{0.0}\OperatorTok{,} \FloatTok{100.0}\NormalTok{)}\OperatorTok{,}\NormalTok{ (}\FloatTok{0.0}\OperatorTok{,} \FloatTok{100.0}\NormalTok{)}\OperatorTok{,}\NormalTok{ (}\FloatTok{0.0}\OperatorTok{,} \FloatTok{100.0}\NormalTok{)}\OperatorTok{,}\NormalTok{ (}\FloatTok{0.0}\OperatorTok{,} \FloatTok{120.0}\NormalTok{))}
\NormalTok{result }\OperatorTok{=}\NormalTok{ classify\_system(ds}\OperatorTok{,}\NormalTok{ bounds)}
\NormalTok{print\_classification(result)}
\end{Highlighting}
\end{Shaded}

\hypertarget{related-software}{%
\section{Related Software}\label{related-software}}

Several Julia packages provide complementary functionality.
\texttt{Differential\-Equations.jl} \cite{Rackauckas2017} offers
comprehensive ODE/SDE solvers but does not perform structural
classification. \texttt{BifurcationKit.jl} focuses on continuation and
bifurcation analysis. \texttt{DynamicalSystems.jl} \cite{Datseris2018}
provides tools for chaos detection and attractor characterisation.
\texttt{GAIO.jl} \cite{GAIO2025} provides numerical routines for the
global analysis of dynamical systems. \texttt{FlowClass.jl} complements
these by addressing a distinct question: not \emph{how} a system
behaves, but \emph{what kind} of system it is structurally.

In Python, \texttt{PyDSTool} and \texttt{PySCeS} offer dynamical systems
modelling for biology, but neither provides systematic structural
classification. The landscape and flux decomposition methods of
@Wang2015 are related theoretically; \texttt{FlowClass.jl} provides
practical tools for the classification aspect of this framework.

\hypertarget{conclusion}{%
\section{Conclusion}\label{conclusion}}

Our ability to make qualitative statements about a dynamical system
depends crucially on the nature of the dynamics. The qualitative aspects
can deliver profound biological systems. If stem cell differentiation
were to follow gradient dynamics, for example, then the forward and
backward paths through gene expression space would be identical
\cite{Guillemin2021,Vittadello2025}. For gradient systems, and
gradient-like systems we have access to Lyapunov functions, and concepts
from catastrophe theory can be applied and yield powerful insights
\cite{Rand2021}.

For more general dynamical systems we cannot fall back on such elegant
theory. A main focus in the design of\texttt{FlowClass.jl} was to
classify dynamical systems into the relevant categories that determine
whether or not a given system (at least in the specific parameterisation
considered). The routines provided as part of the package make

\hypertarget{acknowledgements}{%
\section{Acknowledgements}\label{acknowledgements}}

The author thanks the Australian Research Council for financial support
through an ARC Laureate Fellowship (FL220100005) \# References

\bibliographystyle{unsrt}
\bibliography{paper}

\appendix
\newpage
{\LARGE \bf APPENDIX}
\par

The Appendix contains the README file of the \texttt{FlowClass.jl} package. The package, including all examples can be downloaded from \url{https://github.com/theosysbio/FlowClass.jl}.

\hypertarget{flowclass.jl}{%
\section{FlowClass.jl}\label{flowclass.jl}}

A Julia package for classifying continuous-time dynamical systems into a
hierarchy of structural classes: Gradient, Gradient-like, Morse-Smale,
Structurally Stable, and General.

\hypertarget{motivation}{%
\subsection{Motivation}\label{motivation}}

Many biological and physical systems are modelled as dynamical systems
of the form:

\[\frac{d\mathbf{x}}{dt} = \mathbf{F}(\mathbf{x})\]

Understanding which structural class a system belongs to has profound
implications for its qualitative behaviour. For instance:

\begin{itemize}
\tightlist
\item
  \textbf{Gradient systems} cannot exhibit oscillations or chaos ---
  trajectories always descend a potential
\item
  \textbf{Morse-Smale systems} can have limit cycles but remain
  structurally stable in less than 3 dimensions
\item
  \textbf{General systems} may exhibit chaotic dynamics or complex
  attractors
\end{itemize}

This package provides computational tools to test and classify dynamical
systems based on their structural properties, with applications to
Waddington's epigenetic landscape and cell fate decision models.

\hypertarget{installation}{%
\subsection{Installation}\label{installation}}

\begin{Shaded}
\begin{Highlighting}[]
\KeywordTok{using}\NormalTok{ Pkg}
\NormalTok{Pkg.add(url}\OperatorTok{=}\StringTok{"https://github.com/Theosysbio/FlowClass.jl"}\NormalTok{)}
\end{Highlighting}
\end{Shaded}

Or for local development:

\begin{Shaded}
\begin{Highlighting}[]
\NormalTok{cd(}\StringTok{"path/to/FlowClass.jl"}\NormalTok{)}
\KeywordTok{using}\NormalTok{ Pkg}
\NormalTok{Pkg.activate(}\StringTok{"."}\NormalTok{)}
\NormalTok{Pkg.instantiate()}
\end{Highlighting}
\end{Shaded}

\hypertarget{quick-start}{%
\subsection{Quick Start}\label{quick-start}}

\begin{Shaded}
\begin{Highlighting}[]
\KeywordTok{using}\NormalTok{ FlowClass}

\CommentTok{\# Define a dynamical system: dx/dt = F(x)}
\CommentTok{\# Example: A gradient system with potential V(x) = x₁² + x₁x₂ + x₂²}
\NormalTok{F }\OperatorTok{=}\NormalTok{ x }\OperatorTok{{-}\textgreater{}}\NormalTok{ [}\OperatorTok{{-}}\FloatTok{2}\NormalTok{x[}\FloatTok{1}\NormalTok{] }\OperatorTok{{-}}\NormalTok{ x[}\FloatTok{2}\NormalTok{]}\OperatorTok{,} \OperatorTok{{-}}\NormalTok{x[}\FloatTok{1}\NormalTok{] }\OperatorTok{{-}} \FloatTok{2}\NormalTok{x[}\FloatTok{2}\NormalTok{]]}
\NormalTok{ds }\OperatorTok{=}\NormalTok{ DynamicalSystem(F}\OperatorTok{,} \FloatTok{2}\NormalTok{)}

\CommentTok{\# Classify the system}
\NormalTok{bounds }\OperatorTok{=}\NormalTok{ ((}\OperatorTok{{-}}\FloatTok{2.0}\OperatorTok{,} \FloatTok{2.0}\NormalTok{)}\OperatorTok{,}\NormalTok{ (}\OperatorTok{{-}}\FloatTok{2.0}\OperatorTok{,} \FloatTok{2.0}\NormalTok{))}
\NormalTok{result }\OperatorTok{=}\NormalTok{ classify\_system(ds}\OperatorTok{,}\NormalTok{ bounds)}
\NormalTok{print\_classification(result)}
\end{Highlighting}
\end{Shaded}

\hypertarget{the-classification-hierarchy}{%
\subsection{The Classification
Hierarchy}\label{the-classification-hierarchy}}

From most restrictive to most general:

\begin{verbatim}
Gradient Systems
└── Gradient-like Systems
    └── Morse-Smale Systems
        └── Structurally Stable Systems
            └── General Dynamical Systems
\end{verbatim}

\hypertarget{key-properties-by-class}{%
\subsubsection{Key Properties by Class}\label{key-properties-by-class}}

\begin{longtable}[]{@{}llllll@{}}
\toprule
Class & Jacobian & Curl & Periodic Orbits & Fixed Points & Lyapunov
Function \\
\midrule
\endhead
Gradient & Symmetric & Zero & None & All hyperbolic & Global \\
Gradient-like & Nearly symmetric & Near zero & None & All hyperbolic &
Global \\
Morse-Smale & No requirement & Any & Hyperbolic only & All hyperbolic &
Local only \\
Structurally Stable & No requirement & Any & Hyperbolic only & All
hyperbolic & Local \\
General & No requirement & Any & Any & Non-hyperbolic possible & None
guaranteed \\
\bottomrule
\end{longtable}

\hypertarget{systemclass-enum}{%
\subsubsection{SystemClass Enum}\label{systemclass-enum}}

\begin{Shaded}
\begin{Highlighting}[]
\PreprocessorTok{@enum}\NormalTok{ SystemClass }\KeywordTok{begin}
\NormalTok{    GRADIENT           }\CommentTok{\# Pure gradient system: F = {-}∇V}
\NormalTok{    GRADIENT\_LIKE      }\CommentTok{\# Has global Lyapunov function, no periodic orbits}
\NormalTok{    MORSE\_SMALE        }\CommentTok{\# Hyperbolic fixed points and orbits, transverse manifolds}
\NormalTok{    STRUCTURALLY\_STABLE }\CommentTok{\# Robust to small perturbations}
\NormalTok{    GENERAL            }\CommentTok{\# No special structure guaranteed}
\KeywordTok{end}
\end{Highlighting}
\end{Shaded}

\hypertarget{api-reference}{%
\subsection{API Reference}\label{api-reference}}

\hypertarget{types}{%
\subsubsection{Types}\label{types}}

\hypertarget{dynamicalsystemf}{%
\paragraph{\texorpdfstring{\texttt{DynamicalSystem\{F\}}}{DynamicalSystem\{F\}}}\label{dynamicalsystemf}}

Represents a continuous-time dynamical system dx/dt = f(x).

\begin{Shaded}
\begin{Highlighting}[]
\CommentTok{\# From function and dimension}
\NormalTok{ds }\OperatorTok{=}\NormalTok{ DynamicalSystem(x }\OperatorTok{{-}\textgreater{}} \OperatorTok{{-}}\NormalTok{x}\OperatorTok{,} \FloatTok{2}\NormalTok{)}

\CommentTok{\# From function and sample point (infers dimension)}
\NormalTok{ds }\OperatorTok{=}\NormalTok{ DynamicalSystem(x }\OperatorTok{{-}\textgreater{}} \OperatorTok{{-}}\NormalTok{x}\OperatorTok{,}\NormalTok{ [}\FloatTok{1.0}\OperatorTok{,} \FloatTok{2.0}\NormalTok{])}

\CommentTok{\# Evaluate the vector field}
\NormalTok{ds([}\FloatTok{1.0}\OperatorTok{,} \FloatTok{2.0}\NormalTok{])  }\CommentTok{\# returns [{-}1.0, {-}2.0]}

\CommentTok{\# Get dimension}
\NormalTok{dimension(ds)  }\CommentTok{\# returns 2}
\end{Highlighting}
\end{Shaded}

\hypertarget{fixedpoint}{%
\paragraph{\texorpdfstring{\texttt{FixedPoint}}{FixedPoint}}\label{fixedpoint}}

Represents a fixed point with stability information.

\begin{Shaded}
\begin{Highlighting}[]
\KeywordTok{struct}\NormalTok{ FixedPoint}
\NormalTok{    location}\OperatorTok{::}\DataTypeTok{Vector}\NormalTok{\{}\DataTypeTok{Float64}\NormalTok{\}    }\CommentTok{\# Position in state space}
\NormalTok{    eigenvalues}\OperatorTok{::}\DataTypeTok{Vector}\NormalTok{\{}\DataTypeTok{ComplexF64}\NormalTok{\}  }\CommentTok{\# Eigenvalues of Jacobian}
    \KeywordTok{type}\OperatorTok{::}\NormalTok{FixedPointType         }\CommentTok{\# Classification}
\KeywordTok{end}
\end{Highlighting}
\end{Shaded}

\hypertarget{fixedpointtype-enum}{%
\paragraph{\texorpdfstring{\texttt{FixedPointType}
Enum}{FixedPointType Enum}}\label{fixedpointtype-enum}}

\begin{Shaded}
\begin{Highlighting}[]
\PreprocessorTok{@enum}\NormalTok{ FixedPointType }\KeywordTok{begin}
\NormalTok{    STABLE\_NODE        }\CommentTok{\# All eigenvalues have negative real parts (no imaginary)}
\NormalTok{    UNSTABLE\_NODE      }\CommentTok{\# All eigenvalues have positive real parts (no imaginary)}
\NormalTok{    SADDLE             }\CommentTok{\# Mixed signs of real parts}
\NormalTok{    STABLE\_FOCUS       }\CommentTok{\# Negative real parts with imaginary components}
\NormalTok{    UNSTABLE\_FOCUS     }\CommentTok{\# Positive real parts with imaginary components}
\NormalTok{    CENTER             }\CommentTok{\# Pure imaginary eigenvalues}
\NormalTok{    NON\_HYPERBOLIC     }\CommentTok{\# At least one eigenvalue with zero real part}
\KeywordTok{end}
\end{Highlighting}
\end{Shaded}

\hypertarget{periodicorbit}{%
\paragraph{\texorpdfstring{\texttt{PeriodicOrbit}}{PeriodicOrbit}}\label{periodicorbit}}

Represents a detected periodic orbit.

\begin{Shaded}
\begin{Highlighting}[]
\KeywordTok{struct}\NormalTok{ PeriodicOrbit}
\NormalTok{    points}\OperatorTok{::}\DataTypeTok{Vector}\NormalTok{\{}\DataTypeTok{Vector}\NormalTok{\{}\DataTypeTok{Float64}\NormalTok{\}\}  }\CommentTok{\# Sample points along orbit}
\NormalTok{    period}\OperatorTok{::}\DataTypeTok{Float64}                   \CommentTok{\# Estimated period}
\NormalTok{    is\_stable}\OperatorTok{::}\DataTypeTok{Bool}                   \CommentTok{\# Stability (via Floquet analysis)}
\KeywordTok{end}
\end{Highlighting}
\end{Shaded}

\hypertarget{classificationresult}{%
\paragraph{\texorpdfstring{\texttt{ClassificationResult}}{ClassificationResult}}\label{classificationresult}}

Complete result from system classification.

\begin{Shaded}
\begin{Highlighting}[]
\KeywordTok{struct}\NormalTok{ ClassificationResult}
\NormalTok{    system\_class}\OperatorTok{::}\NormalTok{SystemClass}
\NormalTok{    fixed\_points}\OperatorTok{::}\DataTypeTok{Vector}\NormalTok{\{FixedPoint\}}
\NormalTok{    periodic\_orbits}\OperatorTok{::}\DataTypeTok{Vector}\NormalTok{\{PeriodicOrbit\}}
\NormalTok{    jacobian\_symmetry}\OperatorTok{::}\DataTypeTok{Float64}       \CommentTok{\# Mean relative symmetry error}
\NormalTok{    curl\_gradient\_ratio}\OperatorTok{::}\DataTypeTok{Float64}     \CommentTok{\# ‖curl‖ / ‖gradient‖}
\NormalTok{    has\_transverse\_manifolds}\OperatorTok{::}\DataTypeTok{Union}\NormalTok{\{}\DataTypeTok{Bool}\OperatorTok{,} \DataTypeTok{Nothing}\NormalTok{\}}
\NormalTok{    confidence}\OperatorTok{::}\DataTypeTok{Float64}              \CommentTok{\# Classification confidence}
\NormalTok{    details}\OperatorTok{::}\DataTypeTok{Dict}\NormalTok{\{}\DataTypeTok{String}\OperatorTok{,} \DataTypeTok{Any}\NormalTok{\}       }\CommentTok{\# Additional analysis data}
\KeywordTok{end}
\end{Highlighting}
\end{Shaded}

\hypertarget{jacobian-analysis}{%
\subsubsection{Jacobian Analysis}\label{jacobian-analysis}}

\hypertarget{compute_jacobiands-x-compute_jacobianf-x}{%
\paragraph{\texorpdfstring{\texttt{compute\_jacobian(ds,\ x)} /
\texttt{compute\_jacobian(f,\ x)}}{compute\_jacobian(ds, x) / compute\_jacobian(f, x)}}\label{compute_jacobiands-x-compute_jacobianf-x}}

Compute the Jacobian matrix J{[}i,j{]} = ∂fᵢ/∂xⱼ at point x using
automatic differentiation.

\begin{Shaded}
\begin{Highlighting}[]
\NormalTok{ds }\OperatorTok{=}\NormalTok{ DynamicalSystem(x }\OperatorTok{{-}\textgreater{}}\NormalTok{ [x[}\FloatTok{1}\NormalTok{]}\OperatorTok{\^{}}\FloatTok{2}\OperatorTok{,}\NormalTok{ x[}\FloatTok{1}\NormalTok{]}\OperatorTok{*}\NormalTok{x[}\FloatTok{2}\NormalTok{]]}\OperatorTok{,} \FloatTok{2}\NormalTok{)}
\NormalTok{J }\OperatorTok{=}\NormalTok{ compute\_jacobian(ds}\OperatorTok{,}\NormalTok{ [}\FloatTok{2.0}\OperatorTok{,} \FloatTok{3.0}\NormalTok{])}
\CommentTok{\# J = [4.0 0.0; 3.0 2.0]}
\end{Highlighting}
\end{Shaded}

\hypertarget{is_jacobian_symmetricj-rtol1e-8-atol1e-10}{%
\paragraph{\texorpdfstring{\texttt{is\_jacobian\_symmetric(J;\ rtol=1e-8,\ atol=1e-10)}}{is\_jacobian\_symmetric(J; rtol=1e-8, atol=1e-10)}}\label{is_jacobian_symmetricj-rtol1e-8-atol1e-10}}

Test whether a Jacobian matrix is symmetric within tolerance.

\begin{Shaded}
\begin{Highlighting}[]
\NormalTok{J\_sym }\OperatorTok{=}\NormalTok{ [}\OperatorTok{{-}}\FloatTok{2.0} \FloatTok{0.5}\OperatorTok{;} \FloatTok{0.5} \OperatorTok{{-}}\FloatTok{1.0}\NormalTok{]}
\NormalTok{is\_jacobian\_symmetric(J\_sym)  }\CommentTok{\# true}

\NormalTok{J\_nonsym }\OperatorTok{=}\NormalTok{ [}\OperatorTok{{-}}\FloatTok{1.0} \FloatTok{0.5}\OperatorTok{;} \OperatorTok{{-}}\FloatTok{0.5} \OperatorTok{{-}}\FloatTok{1.0}\NormalTok{]}
\NormalTok{is\_jacobian\_symmetric(J\_nonsym)  }\CommentTok{\# false}
\end{Highlighting}
\end{Shaded}

\hypertarget{jacobian_symmetry_errorj-jacobian_symmetry_errords-x}{%
\paragraph{\texorpdfstring{\texttt{jacobian\_symmetry\_error(J)} /
\texttt{jacobian\_symmetry\_error(ds,\ x)}}{jacobian\_symmetry\_error(J) / jacobian\_symmetry\_error(ds, x)}}\label{jacobian_symmetry_errorj-jacobian_symmetry_errords-x}}

Compute the Frobenius norm of the antisymmetric part: ‖(J − Jᵀ)/2‖.

\hypertarget{relative_jacobian_symmetry_errorj-relative_jacobian_symmetry_errords-x}{%
\paragraph{\texorpdfstring{\texttt{relative\_jacobian\_symmetry\_error(J)}
/
\texttt{relative\_jacobian\_symmetry\_error(ds,\ x)}}{relative\_jacobian\_symmetry\_error(J) / relative\_jacobian\_symmetry\_error(ds, x)}}\label{relative_jacobian_symmetry_errorj-relative_jacobian_symmetry_errords-x}}

Scale-independent symmetry error: ‖(J − Jᵀ)/2‖ / ‖J‖.

\hypertarget{curl-analysis}{%
\subsubsection{Curl Analysis}\label{curl-analysis}}

For a vector field \textbf{F}, the curl measures the rotational
component of the dynamics. In the Helmholtz decomposition \textbf{F} =
−∇U + \textbf{F}\_curl, the curl component \textbf{F}\_curl is
orthogonal to the gradient and cannot be captured by any potential
landscape.

\hypertarget{curl_magnitudeds-x-curl_magnitudef-x-n}{%
\paragraph{\texorpdfstring{\texttt{curl\_magnitude(ds,\ x)} /
\texttt{curl\_magnitude(f,\ x,\ n)}}{curl\_magnitude(ds, x) / curl\_magnitude(f, x, n)}}\label{curl_magnitudeds-x-curl_magnitudef-x-n}}

Compute the magnitude of the curl at point x. For 2D systems, returns
the scalar curl. For 3D, returns ‖∇ × \textbf{F}‖. For higher
dimensions, returns ‖(J − Jᵀ)/2‖\_F (Frobenius norm of antisymmetric
part).

\begin{Shaded}
\begin{Highlighting}[]
\CommentTok{\# Rotation system has high curl}
\NormalTok{rotation }\OperatorTok{=}\NormalTok{ DynamicalSystem(x }\OperatorTok{{-}\textgreater{}}\NormalTok{ [}\OperatorTok{{-}}\NormalTok{x[}\FloatTok{2}\NormalTok{]}\OperatorTok{,}\NormalTok{ x[}\FloatTok{1}\NormalTok{]]}\OperatorTok{,} \FloatTok{2}\NormalTok{)}
\NormalTok{curl\_magnitude(rotation}\OperatorTok{,}\NormalTok{ [}\FloatTok{1.0}\OperatorTok{,} \FloatTok{0.0}\NormalTok{])  }\CommentTok{\# ≈ 2.0}

\CommentTok{\# Gradient system has zero curl}
\NormalTok{gradient\_sys }\OperatorTok{=}\NormalTok{ DynamicalSystem(x }\OperatorTok{{-}\textgreater{}}\NormalTok{ [}\OperatorTok{{-}}\FloatTok{2}\NormalTok{x[}\FloatTok{1}\NormalTok{]}\OperatorTok{,} \OperatorTok{{-}}\FloatTok{2}\NormalTok{x[}\FloatTok{2}\NormalTok{]]}\OperatorTok{,} \FloatTok{2}\NormalTok{)}
\NormalTok{curl\_magnitude(gradient\_sys}\OperatorTok{,}\NormalTok{ [}\FloatTok{1.0}\OperatorTok{,} \FloatTok{1.0}\NormalTok{])  }\CommentTok{\# ≈ 0.0}
\end{Highlighting}
\end{Shaded}

\hypertarget{is_curl_freeds-x-atol1e-10-is_curl_freeds-bounds-n_samples100-atol1e-10}{%
\paragraph{\texorpdfstring{\texttt{is\_curl\_free(ds,\ x;\ atol=1e-10)}
/
\texttt{is\_curl\_free(ds,\ bounds;\ n\_samples=100,\ atol=1e-10)}}{is\_curl\_free(ds, x; atol=1e-10) / is\_curl\_free(ds, bounds; n\_samples=100, atol=1e-10)}}\label{is_curl_freeds-x-atol1e-10-is_curl_freeds-bounds-n_samples100-atol1e-10}}

Test if the curl is zero at a point or throughout a region.

\hypertarget{curl_to_gradient_ratiods-x}{%
\paragraph{\texorpdfstring{\texttt{curl\_to\_gradient\_ratio(ds,\ x)}}{curl\_to\_gradient\_ratio(ds, x)}}\label{curl_to_gradient_ratiods-x}}

Compute the ratio ‖curl‖ / ‖F‖, indicating the relative strength of
rotational dynamics.

\hypertarget{fixed-point-analysis}{%
\subsubsection{Fixed Point Analysis}\label{fixed-point-analysis}}

\hypertarget{find_fixed_pointsds-bounds-n_starts100-tol1e-8}{%
\paragraph{\texorpdfstring{\texttt{find\_fixed\_points(ds,\ bounds;\ n\_starts=100,\ tol=1e-8)}}{find\_fixed\_points(ds, bounds; n\_starts=100, tol=1e-8)}}\label{find_fixed_pointsds-bounds-n_starts100-tol1e-8}}

Find fixed points of the system within the specified bounds using
multi-start optimisation.

\begin{Shaded}
\begin{Highlighting}[]
\CommentTok{\# Toggle switch with two stable states}
\NormalTok{toggle }\OperatorTok{=}\NormalTok{ DynamicalSystem(x }\OperatorTok{{-}\textgreater{}}\NormalTok{ [}
    \FloatTok{1}\OperatorTok{/}\NormalTok{(}\FloatTok{1} \OperatorTok{+}\NormalTok{ x[}\FloatTok{2}\NormalTok{]}\OperatorTok{\^{}}\FloatTok{2}\NormalTok{) }\OperatorTok{{-}}\NormalTok{ x[}\FloatTok{1}\NormalTok{]}\OperatorTok{,}
    \FloatTok{1}\OperatorTok{/}\NormalTok{(}\FloatTok{1} \OperatorTok{+}\NormalTok{ x[}\FloatTok{1}\NormalTok{]}\OperatorTok{\^{}}\FloatTok{2}\NormalTok{) }\OperatorTok{{-}}\NormalTok{ x[}\FloatTok{2}\NormalTok{]}
\NormalTok{]}\OperatorTok{,} \FloatTok{2}\NormalTok{)}
\NormalTok{bounds }\OperatorTok{=}\NormalTok{ ((}\FloatTok{0.0}\OperatorTok{,} \FloatTok{2.0}\NormalTok{)}\OperatorTok{,}\NormalTok{ (}\FloatTok{0.0}\OperatorTok{,} \FloatTok{2.0}\NormalTok{))}

\NormalTok{fps }\OperatorTok{=}\NormalTok{ find\_fixed\_points(toggle}\OperatorTok{,}\NormalTok{ bounds)}
\KeywordTok{for}\NormalTok{ fp }\KeywordTok{in}\NormalTok{ fps}
\NormalTok{    println(}\StringTok{"Fixed point at $(fp.location): $(fp.type)"}\NormalTok{)}
\KeywordTok{end}
\end{Highlighting}
\end{Shaded}

\hypertarget{classify_fixed_pointds-x-classify_fixed_pointeigenvalues}{%
\paragraph{\texorpdfstring{\texttt{classify\_fixed\_point(ds,\ x)} /
\texttt{classify\_fixed\_point(eigenvalues)}}{classify\_fixed\_point(ds, x) / classify\_fixed\_point(eigenvalues)}}\label{classify_fixed_pointds-x-classify_fixed_pointeigenvalues}}

Determine the type of a fixed point from its Jacobian eigenvalues.

\hypertarget{is_hyperbolicfpfixedpoint-is_hyperboliceigenvalues}{%
\paragraph{\texorpdfstring{\texttt{is\_hyperbolic(fp::FixedPoint)} /
\texttt{is\_hyperbolic(eigenvalues)}}{is\_hyperbolic(fp::FixedPoint) / is\_hyperbolic(eigenvalues)}}\label{is_hyperbolicfpfixedpoint-is_hyperboliceigenvalues}}

Check if a fixed point is hyperbolic (no eigenvalues with zero real
part).

\hypertarget{periodic-orbit-detection}{%
\subsubsection{Periodic Orbit
Detection}\label{periodic-orbit-detection}}

\hypertarget{find_periodic_orbitsds-bounds-n_trajectories50-max_period100.0}{%
\paragraph{\texorpdfstring{\texttt{find\_periodic\_orbits(ds,\ bounds;\ n\_trajectories=50,\ max\_period=100.0)}}{find\_periodic\_orbits(ds, bounds; n\_trajectories=50, max\_period=100.0)}}\label{find_periodic_orbitsds-bounds-n_trajectories50-max_period100.0}}

Search for periodic orbits by integrating trajectories and detecting
recurrence.

\begin{Shaded}
\begin{Highlighting}[]
\CommentTok{\# Van der Pol oscillator (has a limit cycle)}
\NormalTok{vdp }\OperatorTok{=}\NormalTok{ DynamicalSystem(x }\OperatorTok{{-}\textgreater{}}\NormalTok{ [x[}\FloatTok{2}\NormalTok{]}\OperatorTok{,}\NormalTok{ (}\FloatTok{1} \OperatorTok{{-}}\NormalTok{ x[}\FloatTok{1}\NormalTok{]}\OperatorTok{\^{}}\FloatTok{2}\NormalTok{)}\OperatorTok{*}\NormalTok{x[}\FloatTok{2}\NormalTok{] }\OperatorTok{{-}}\NormalTok{ x[}\FloatTok{1}\NormalTok{]]}\OperatorTok{,} \FloatTok{2}\NormalTok{)}
\NormalTok{bounds }\OperatorTok{=}\NormalTok{ ((}\OperatorTok{{-}}\FloatTok{3.0}\OperatorTok{,} \FloatTok{3.0}\NormalTok{)}\OperatorTok{,}\NormalTok{ (}\OperatorTok{{-}}\FloatTok{3.0}\OperatorTok{,} \FloatTok{3.0}\NormalTok{))}

\NormalTok{orbits }\OperatorTok{=}\NormalTok{ find\_periodic\_orbits(vdp}\OperatorTok{,}\NormalTok{ bounds)}
\KeywordTok{if} \OperatorTok{!}\NormalTok{isempty(orbits)}
\NormalTok{    println(}\StringTok{"Found orbit with period ≈ $(orbits[1].period)"}\NormalTok{)}
\KeywordTok{end}
\end{Highlighting}
\end{Shaded}

\hypertarget{has_periodic_orbitsds-bounds-kwargs...}{%
\paragraph{\texorpdfstring{\texttt{has\_periodic\_orbits(ds,\ bounds;\ kwargs...)}}{has\_periodic\_orbits(ds, bounds; kwargs...)}}\label{has_periodic_orbitsds-bounds-kwargs...}}

Quick check for the existence of periodic orbits. Returns \texttt{true}
if any are found.

\hypertarget{manifold-analysis}{%
\subsubsection{Manifold Analysis}\label{manifold-analysis}}

\hypertarget{compute_stable_manifoldds-fp-n_points100-extent1.0}{%
\paragraph{\texorpdfstring{\texttt{compute\_stable\_manifold(ds,\ fp;\ n\_points=100,\ extent=1.0)}}{compute\_stable\_manifold(ds, fp; n\_points=100, extent=1.0)}}\label{compute_stable_manifoldds-fp-n_points100-extent1.0}}

Compute points along the stable manifold of a saddle point.

\hypertarget{compute_unstable_manifoldds-fp-n_points100-extent1.0}{%
\paragraph{\texorpdfstring{\texttt{compute\_unstable\_manifold(ds,\ fp;\ n\_points=100,\ extent=1.0)}}{compute\_unstable\_manifold(ds, fp; n\_points=100, extent=1.0)}}\label{compute_unstable_manifoldds-fp-n_points100-extent1.0}}

Compute points along the unstable manifold of a saddle point.

\hypertarget{detect_homoclinic_orbitds-saddle-tol0.1}{%
\paragraph{\texorpdfstring{\texttt{detect\_homoclinic\_orbit(ds,\ saddle;\ tol=0.1)}}{detect\_homoclinic\_orbit(ds, saddle; tol=0.1)}}\label{detect_homoclinic_orbitds-saddle-tol0.1}}

Check for homoclinic connections (orbits connecting a saddle to itself).

\hypertarget{check_transversalityds-fps-tol0.01}{%
\paragraph{\texorpdfstring{\texttt{check\_transversality(ds,\ fps;\ tol=0.01)}}{check\_transversality(ds, fps; tol=0.01)}}\label{check_transversalityds-fps-tol0.01}}

Verify that stable and unstable manifolds intersect transversally
(required for Morse-Smale).

\hypertarget{classification-functions}{%
\subsubsection{Classification
Functions}\label{classification-functions}}

\hypertarget{classify_systemds-bounds-kwargs...}{%
\paragraph{\texorpdfstring{\texttt{classify\_system(ds,\ bounds;\ kwargs...)}}{classify\_system(ds, bounds; kwargs...)}}\label{classify_systemds-bounds-kwargs...}}

Perform full classification with detailed analysis.

\begin{Shaded}
\begin{Highlighting}[]
\NormalTok{ds }\OperatorTok{=}\NormalTok{ DynamicalSystem(x }\OperatorTok{{-}\textgreater{}}\NormalTok{ [}\OperatorTok{{-}}\FloatTok{2}\NormalTok{x[}\FloatTok{1}\NormalTok{]}\OperatorTok{,} \OperatorTok{{-}}\FloatTok{3}\NormalTok{x[}\FloatTok{2}\NormalTok{]]}\OperatorTok{,} \FloatTok{2}\NormalTok{)}
\NormalTok{bounds }\OperatorTok{=}\NormalTok{ ((}\OperatorTok{{-}}\FloatTok{2.0}\OperatorTok{,} \FloatTok{2.0}\NormalTok{)}\OperatorTok{,}\NormalTok{ (}\OperatorTok{{-}}\FloatTok{2.0}\OperatorTok{,} \FloatTok{2.0}\NormalTok{))}

\NormalTok{result }\OperatorTok{=}\NormalTok{ classify\_system(ds}\OperatorTok{,}\NormalTok{ bounds)}
\NormalTok{print\_classification(result)}
\end{Highlighting}
\end{Shaded}

\textbf{Keyword arguments:} - \texttt{n\_samples::Int=100} --- Points
sampled for Jacobian/curl analysis - \texttt{n\_starts::Int=100} ---
Starting points for fixed point search -
\texttt{check\_manifolds::Bool=true} --- Whether to analyse manifold
transversality - \texttt{orbit\_timeout::Float64=10.0} --- Max time for
periodic orbit search

\hypertarget{quick_classifyds-bounds}{%
\paragraph{\texorpdfstring{\texttt{quick\_classify(ds,\ bounds)}}{quick\_classify(ds, bounds)}}\label{quick_classifyds-bounds}}

Fast classification with fewer samples and no manifold analysis.

\hypertarget{get_system_classds-bounds}{%
\paragraph{\texorpdfstring{\texttt{get\_system\_class(ds,\ bounds)}}{get\_system\_class(ds, bounds)}}\label{get_system_classds-bounds}}

Return only the \texttt{SystemClass} enum value.

\hypertarget{classification-result-queries}{%
\subsubsection{Classification Result
Queries}\label{classification-result-queries}}

\begin{Shaded}
\begin{Highlighting}[]
\NormalTok{result }\OperatorTok{=}\NormalTok{ classify\_system(ds}\OperatorTok{,}\NormalTok{ bounds)}

\NormalTok{is\_gradient(result)              }\CommentTok{\# true if GRADIENT}
\NormalTok{is\_gradient\_like(result)         }\CommentTok{\# true if GRADIENT or GRADIENT\_LIKE}
\NormalTok{is\_morse\_smale(result)           }\CommentTok{\# true if Morse{-}Smale or more restrictive}
\NormalTok{allows\_periodic\_orbits(result)   }\CommentTok{\# false for gradient{-}like systems}

\CommentTok{\# Get landscape interpretation}
\NormalTok{can\_represent}\OperatorTok{,}\NormalTok{ landscape\_type}\OperatorTok{,}\NormalTok{ description }\OperatorTok{=}\NormalTok{ has\_landscape\_representation(result)}
\end{Highlighting}
\end{Shaded}

\hypertarget{utility-functions}{%
\subsubsection{Utility Functions}\label{utility-functions}}

\hypertarget{print_classificationresult-iostdout}{%
\paragraph{\texorpdfstring{\texttt{print\_classification(result;\ io=stdout)}}{print\_classification(result; io=stdout)}}\label{print_classificationresult-iostdout}}

Print a formatted classification report.

\begin{Shaded}
\begin{Highlighting}[]
\NormalTok{result }\OperatorTok{=}\NormalTok{ classify\_system(ds}\OperatorTok{,}\NormalTok{ bounds)}
\NormalTok{print\_classification(result)}
\end{Highlighting}
\end{Shaded}

Output:

\begin{verbatim}
╔══════════════════════════════════════════════════════════════╗
║                  System Classification Report                 ║
╠══════════════════════════════════════════════════════════════╣
║ System Class: GRADIENT                                       ║
║ Confidence: 0.95                                             ║
╠══════════════════════════════════════════════════════════════╣
║ Fixed Points: 1                                              ║
║   • Stable node at [0.0, 0.0]                               ║
║ Periodic Orbits: 0                                           ║
╠══════════════════════════════════════════════════════════════╣
║ Jacobian Symmetry Error: 1.2e-15                            ║
║ Curl/Gradient Ratio: 0.0                                     ║
║ Manifolds Transverse: N/A (no saddles)                      ║
╠══════════════════════════════════════════════════════════════╣
║ Landscape: Global potential V(x) exists where F = -∇V        ║
╚══════════════════════════════════════════════════════════════╝
\end{verbatim}

\hypertarget{examples}{%
\subsection{Examples}\label{examples}}

\hypertarget{example-1-testing-a-gradient-system}{%
\subsubsection{Example 1: Testing a Gradient
System}\label{example-1-testing-a-gradient-system}}

A gradient system satisfies dx/dt = −∇V(x) for some scalar potential V.
Its Jacobian is the negative Hessian of V, which is always symmetric.

\begin{Shaded}
\begin{Highlighting}[]
\KeywordTok{using}\NormalTok{ FlowClass}

\CommentTok{\# Potential: V(x) = x₁² + x₂² (paraboloid)}
\CommentTok{\# Gradient: ∇V = [2x₁, 2x₂]}
\CommentTok{\# System: dx/dt = {-}∇V = [{-}2x₁, {-}2x₂]}

\NormalTok{ds }\OperatorTok{=}\NormalTok{ DynamicalSystem(x }\OperatorTok{{-}\textgreater{}} \OperatorTok{{-}}\FloatTok{2} \OperatorTok{.*}\NormalTok{ x}\OperatorTok{,} \FloatTok{2}\NormalTok{)}
\NormalTok{bounds }\OperatorTok{=}\NormalTok{ ((}\OperatorTok{{-}}\FloatTok{2.0}\OperatorTok{,} \FloatTok{2.0}\NormalTok{)}\OperatorTok{,}\NormalTok{ (}\OperatorTok{{-}}\FloatTok{2.0}\OperatorTok{,} \FloatTok{2.0}\NormalTok{))}

\NormalTok{result }\OperatorTok{=}\NormalTok{ classify\_system(ds}\OperatorTok{,}\NormalTok{ bounds)}
\NormalTok{println(}\StringTok{"Class: "}\OperatorTok{,}\NormalTok{ result.system\_class)  }\CommentTok{\# GRADIENT}
\NormalTok{println(}\StringTok{"Symmetry error: "}\OperatorTok{,}\NormalTok{ result.jacobian\_symmetry)  }\CommentTok{\# ≈ 0}
\NormalTok{println(}\StringTok{"Curl ratio: "}\OperatorTok{,}\NormalTok{ result.curl\_gradient\_ratio)  }\CommentTok{\# ≈ 0}
\end{Highlighting}
\end{Shaded}

\hypertarget{example-2-system-with-rotation-non-gradient}{%
\subsubsection{Example 2: System with Rotation
(Non-Gradient)}\label{example-2-system-with-rotation-non-gradient}}

Systems with rotational dynamics have antisymmetric components in their
Jacobian and non-zero curl.

\begin{Shaded}
\begin{Highlighting}[]
\KeywordTok{using}\NormalTok{ FlowClass}

\CommentTok{\# Damped oscillator with rotation}
\CommentTok{\# dx₁/dt = {-}x₁ + ωx₂}
\CommentTok{\# dx₂/dt = {-}ωx₁ {-} x₂}
\NormalTok{ω }\OperatorTok{=} \FloatTok{1.0}
\NormalTok{ds }\OperatorTok{=}\NormalTok{ DynamicalSystem(x }\OperatorTok{{-}\textgreater{}}\NormalTok{ [}\OperatorTok{{-}}\NormalTok{x[}\FloatTok{1}\NormalTok{] }\OperatorTok{+}\NormalTok{ ω}\OperatorTok{*}\NormalTok{x[}\FloatTok{2}\NormalTok{]}\OperatorTok{,} \OperatorTok{{-}}\NormalTok{ω}\OperatorTok{*}\NormalTok{x[}\FloatTok{1}\NormalTok{] }\OperatorTok{{-}}\NormalTok{ x[}\FloatTok{2}\NormalTok{]]}\OperatorTok{,} \FloatTok{2}\NormalTok{)}

\NormalTok{J }\OperatorTok{=}\NormalTok{ compute\_jacobian(ds}\OperatorTok{,}\NormalTok{ [}\FloatTok{0.0}\OperatorTok{,} \FloatTok{0.0}\NormalTok{])}
\CommentTok{\# J = [{-}1  1; {-}1  {-}1]}

\NormalTok{is\_jacobian\_symmetric(J)  }\CommentTok{\# false}
\NormalTok{relative\_jacobian\_symmetry\_error(J)  }\CommentTok{\# ≈ 0.5}
\NormalTok{curl\_magnitude(ds}\OperatorTok{,}\NormalTok{ [}\FloatTok{0.0}\OperatorTok{,} \FloatTok{0.0}\NormalTok{])  }\CommentTok{\# ≈ 2.0}
\end{Highlighting}
\end{Shaded}

\hypertarget{example-3-lorenz-system}{%
\subsubsection{Example 3: Lorenz System}\label{example-3-lorenz-system}}

The Lorenz system is a classic example of a chaotic, non-gradient
system.

\begin{Shaded}
\begin{Highlighting}[]
\KeywordTok{using}\NormalTok{ FlowClass}

\KeywordTok{function}\NormalTok{ lorenz(x}\OperatorTok{;}\NormalTok{ σ}\OperatorTok{=}\FloatTok{10.0}\OperatorTok{,}\NormalTok{ ρ}\OperatorTok{=}\FloatTok{28.0}\OperatorTok{,}\NormalTok{ β}\OperatorTok{=}\FloatTok{8}\OperatorTok{/}\FloatTok{3}\NormalTok{)}
    \KeywordTok{return}\NormalTok{ [σ }\OperatorTok{*}\NormalTok{ (x[}\FloatTok{2}\NormalTok{] }\OperatorTok{{-}}\NormalTok{ x[}\FloatTok{1}\NormalTok{])}\OperatorTok{,}
\NormalTok{            x[}\FloatTok{1}\NormalTok{] }\OperatorTok{*}\NormalTok{ (ρ }\OperatorTok{{-}}\NormalTok{ x[}\FloatTok{3}\NormalTok{]) }\OperatorTok{{-}}\NormalTok{ x[}\FloatTok{2}\NormalTok{]}\OperatorTok{,}
\NormalTok{            x[}\FloatTok{1}\NormalTok{] }\OperatorTok{*}\NormalTok{ x[}\FloatTok{2}\NormalTok{] }\OperatorTok{{-}}\NormalTok{ β }\OperatorTok{*}\NormalTok{ x[}\FloatTok{3}\NormalTok{]]}
\KeywordTok{end}

\NormalTok{ds }\OperatorTok{=}\NormalTok{ DynamicalSystem(lorenz}\OperatorTok{,} \FloatTok{3}\NormalTok{)}
\NormalTok{bounds }\OperatorTok{=}\NormalTok{ ((}\OperatorTok{{-}}\FloatTok{20.0}\OperatorTok{,} \FloatTok{20.0}\NormalTok{)}\OperatorTok{,}\NormalTok{ (}\OperatorTok{{-}}\FloatTok{30.0}\OperatorTok{,} \FloatTok{30.0}\NormalTok{)}\OperatorTok{,}\NormalTok{ (}\FloatTok{0.0}\OperatorTok{,} \FloatTok{50.0}\NormalTok{))}

\NormalTok{result }\OperatorTok{=}\NormalTok{ classify\_system(ds}\OperatorTok{,}\NormalTok{ bounds)}
\NormalTok{println(}\StringTok{"Class: "}\OperatorTok{,}\NormalTok{ result.system\_class)  }\CommentTok{\# GENERAL}
\NormalTok{println(}\StringTok{"Fixed points found: "}\OperatorTok{,}\NormalTok{ length(result.fixed\_points))}
\end{Highlighting}
\end{Shaded}

\hypertarget{example-4-stem-cell-differentiation-model}{%
\subsubsection{Example 4: Stem Cell Differentiation
Model}\label{example-4-stem-cell-differentiation-model}}

This example implements the stem cell differentiation model from
Brackston, Lakatos \& Stumpf (2018), which describes the dynamics of
pluripotency factors (Nanog, Oct4-Sox2, Fgf4) and differentiation marker
(Gata6) under the influence of LIF signalling.

The model demonstrates non-gradient dynamics with curl components,
multiple stable states (pluripotent and differentiated), and transition
states --- key features of Waddington's epigenetic landscape.

\hypertarget{the-model-equations-eqns.-816}{%
\paragraph{The Model Equations (Eqns.
8--16)}\label{the-model-equations-eqns.-816}}

The developmental model consists of four molecular species: Nanog
(\(N\)), Oct4-Sox2 complex (\(O\)), Fgf4 (\(F\)), and Gata6 (\(G\)),
with LIF (\(L\)) as an external control parameter. Under a
quasi-equilibrium assumption, the dynamics are governed by eight
reactions: four production propensities and four degradation
propensities.

\textbf{Production propensities:}

\[a_1 = \frac{k_0 O (k_1 + k_2 N^2 + k_0 O + k_3 L)}{1 + k_0 O (k_2 N^2 + k_0 O + k_3 L + k_4 F^2) + k_5 O G^2} \tag{8}\]

\[a_2 = \frac{k_6 + k_7 O}{1 + k_7 O + k_8 G^2} \tag{9}\]

\[a_3 = \frac{k_9 + k_{10} O}{1 + k_{10} O} \tag{10}\]

\[a_4 = \frac{k_{11} + k_{12} G^2 + k_{14} O}{1 + k_{12} G^2 + k_{13} N^2 + k_{14} O} \tag{11}\]

\textbf{Degradation propensities} (first-order with rate \(k_d\)):

\[a_5 = k_d N \tag{12}\]

\[a_6 = k_d O \tag{13}\]

\[a_7 = k_d F \tag{14}\]

\[a_8 = k_d G \tag{15}\]

\textbf{Stoichiometry matrix:}

The system evolution is described by dx/dt = S · a(x), where the
stoichiometry matrix is:

\[S = \begin{bmatrix} 1 & 0 & 0 & 0 & -1 & 0 & 0 & 0 \\ 0 & 1 & 0 & 0 & 0 & -1 & 0 & 0 \\ 0 & 0 & 1 & 0 & 0 & 0 & -1 & 0 \\ 0 & 0 & 0 & 1 & 0 & 0 & 0 & -1 \end{bmatrix} \tag{16}\]

This yields the ODEs: \(\dot{N} = a_1 - k_d N\),
\(\dot{O} = a_2 - k_d O\), \(\dot{F} = a_3 - k_d F\),
\(\dot{G} = a_4 - k_d G\).

\begin{Shaded}
\begin{Highlighting}[]
\KeywordTok{using}\NormalTok{ FlowClass}

\CommentTok{\# Parameters from Brackston et al. (2018) Table in Methods section}
\KeywordTok{const}\NormalTok{ k }\OperatorTok{=}\NormalTok{ (}
\NormalTok{    k0 }\OperatorTok{=} \FloatTok{0.005}\OperatorTok{,}\NormalTok{ k1 }\OperatorTok{=} \FloatTok{0.01}\OperatorTok{,}\NormalTok{ k2 }\OperatorTok{=} \FloatTok{0.4}\OperatorTok{,}\NormalTok{ k3 }\OperatorTok{=} \FloatTok{1.0}\OperatorTok{,}\NormalTok{ k4 }\OperatorTok{=} \FloatTok{0.1}\OperatorTok{,}
\NormalTok{    k5 }\OperatorTok{=} \FloatTok{0.00135}\OperatorTok{,}\NormalTok{ k6 }\OperatorTok{=} \FloatTok{0.01}\OperatorTok{,}\NormalTok{ k7 }\OperatorTok{=} \FloatTok{0.01}\OperatorTok{,}\NormalTok{ k8 }\OperatorTok{=} \FloatTok{1.0}\OperatorTok{,}\NormalTok{ k9 }\OperatorTok{=} \FloatTok{1.0}\OperatorTok{,}
\NormalTok{    k10 }\OperatorTok{=} \FloatTok{0.01}\OperatorTok{,}\NormalTok{ k11 }\OperatorTok{=} \FloatTok{5.0}\OperatorTok{,}\NormalTok{ k12 }\OperatorTok{=} \FloatTok{1.0}\OperatorTok{,}\NormalTok{ k13 }\OperatorTok{=} \FloatTok{0.005}\OperatorTok{,}\NormalTok{ k14 }\OperatorTok{=} \FloatTok{1.0}\OperatorTok{,}
\NormalTok{    kd }\OperatorTok{=} \FloatTok{1.0}
\NormalTok{)}

\StringTok{"""}
\StringTok{Stem cell differentiation model (Brackston et al. 2018, Eqns. 8–16).}
\StringTok{State vector: x = [N, O, F, G] where}
\StringTok{  N = Nanog, O = Oct4{-}Sox2, F = Fgf4, G = Gata6}
\StringTok{Parameter L controls LIF concentration (external signal).}
\StringTok{"""}
\KeywordTok{function}\NormalTok{ stem\_cell\_model(x}\OperatorTok{;}\NormalTok{ L}\OperatorTok{=}\FloatTok{50.0}\OperatorTok{,}\NormalTok{ p}\OperatorTok{=}\NormalTok{k)}
\NormalTok{    N}\OperatorTok{,}\NormalTok{ O}\OperatorTok{,}\NormalTok{ F}\OperatorTok{,}\NormalTok{ G }\OperatorTok{=}\NormalTok{ x}
    
    \CommentTok{\# Production propensities (Eqns. 8–11)}
\NormalTok{    a1 }\OperatorTok{=}\NormalTok{ p.k0 }\OperatorTok{*}\NormalTok{ O }\OperatorTok{*}\NormalTok{ (p.k1 }\OperatorTok{+}\NormalTok{ p.k2}\OperatorTok{*}\NormalTok{N}\OperatorTok{\^{}}\FloatTok{2} \OperatorTok{+}\NormalTok{ p.k0}\OperatorTok{*}\NormalTok{O }\OperatorTok{+}\NormalTok{ p.k3}\OperatorTok{*}\NormalTok{L) }\OperatorTok{/} 
\NormalTok{         (}\FloatTok{1} \OperatorTok{+}\NormalTok{ p.k0}\OperatorTok{*}\NormalTok{O}\OperatorTok{*}\NormalTok{(p.k2}\OperatorTok{*}\NormalTok{N}\OperatorTok{\^{}}\FloatTok{2} \OperatorTok{+}\NormalTok{ p.k0}\OperatorTok{*}\NormalTok{O }\OperatorTok{+}\NormalTok{ p.k3}\OperatorTok{*}\NormalTok{L }\OperatorTok{+}\NormalTok{ p.k4}\OperatorTok{*}\NormalTok{F}\OperatorTok{\^{}}\FloatTok{2}\NormalTok{) }\OperatorTok{+}\NormalTok{ p.k5}\OperatorTok{*}\NormalTok{O}\OperatorTok{*}\NormalTok{G}\OperatorTok{\^{}}\FloatTok{2}\NormalTok{)}
    
\NormalTok{    a2 }\OperatorTok{=}\NormalTok{ (p.k6 }\OperatorTok{+}\NormalTok{ p.k7}\OperatorTok{*}\NormalTok{O) }\OperatorTok{/}\NormalTok{ (}\FloatTok{1} \OperatorTok{+}\NormalTok{ p.k7}\OperatorTok{*}\NormalTok{O }\OperatorTok{+}\NormalTok{ p.k8}\OperatorTok{*}\NormalTok{G}\OperatorTok{\^{}}\FloatTok{2}\NormalTok{)}
    
\NormalTok{    a3 }\OperatorTok{=}\NormalTok{ (p.k9 }\OperatorTok{+}\NormalTok{ p.k10}\OperatorTok{*}\NormalTok{O) }\OperatorTok{/}\NormalTok{ (}\FloatTok{1} \OperatorTok{+}\NormalTok{ p.k10}\OperatorTok{*}\NormalTok{O)}
    
\NormalTok{    a4 }\OperatorTok{=}\NormalTok{ (p.k11 }\OperatorTok{+}\NormalTok{ p.k12}\OperatorTok{*}\NormalTok{G}\OperatorTok{\^{}}\FloatTok{2} \OperatorTok{+}\NormalTok{ p.k14}\OperatorTok{*}\NormalTok{O) }\OperatorTok{/}\NormalTok{ (}\FloatTok{1} \OperatorTok{+}\NormalTok{ p.k12}\OperatorTok{*}\NormalTok{G}\OperatorTok{\^{}}\FloatTok{2} \OperatorTok{+}\NormalTok{ p.k13}\OperatorTok{*}\NormalTok{N}\OperatorTok{\^{}}\FloatTok{2} \OperatorTok{+}\NormalTok{ p.k14}\OperatorTok{*}\NormalTok{O)}
    
    \CommentTok{\# Net rates: production − degradation (from stoichiometry, Eq. 16)}
\NormalTok{    dN }\OperatorTok{=}\NormalTok{ a1 }\OperatorTok{{-}}\NormalTok{ p.kd }\OperatorTok{*}\NormalTok{ N}
\NormalTok{    dO }\OperatorTok{=}\NormalTok{ a2 }\OperatorTok{{-}}\NormalTok{ p.kd }\OperatorTok{*}\NormalTok{ O}
\NormalTok{    dF }\OperatorTok{=}\NormalTok{ a3 }\OperatorTok{{-}}\NormalTok{ p.kd }\OperatorTok{*}\NormalTok{ F}
\NormalTok{    dG }\OperatorTok{=}\NormalTok{ a4 }\OperatorTok{{-}}\NormalTok{ p.kd }\OperatorTok{*}\NormalTok{ G}
    
    \KeywordTok{return}\NormalTok{ [dN}\OperatorTok{,}\NormalTok{ dO}\OperatorTok{,}\NormalTok{ dF}\OperatorTok{,}\NormalTok{ dG]}
\KeywordTok{end}

\CommentTok{\# Create system with high LIF (favours pluripotency)}
\NormalTok{ds\_high\_LIF }\OperatorTok{=}\NormalTok{ DynamicalSystem(x }\OperatorTok{{-}\textgreater{}}\NormalTok{ stem\_cell\_model(x}\OperatorTok{;}\NormalTok{ L}\OperatorTok{=}\FloatTok{150.0}\NormalTok{)}\OperatorTok{,} \FloatTok{4}\NormalTok{)}

\CommentTok{\# Create system with low LIF (favours differentiation)}
\NormalTok{ds\_low\_LIF }\OperatorTok{=}\NormalTok{ DynamicalSystem(x }\OperatorTok{{-}\textgreater{}}\NormalTok{ stem\_cell\_model(x}\OperatorTok{;}\NormalTok{ L}\OperatorTok{=}\FloatTok{10.0}\NormalTok{)}\OperatorTok{,} \FloatTok{4}\NormalTok{)}

\CommentTok{\# Define bounds for the four{-}dimensional state space}
\CommentTok{\# N ∈ [0, 100], O ∈ [0, 100], F ∈ [0, 100], G ∈ [0, 120]}
\NormalTok{bounds }\OperatorTok{=}\NormalTok{ ((}\FloatTok{0.0}\OperatorTok{,} \FloatTok{100.0}\NormalTok{)}\OperatorTok{,}\NormalTok{ (}\FloatTok{0.0}\OperatorTok{,} \FloatTok{100.0}\NormalTok{)}\OperatorTok{,}\NormalTok{ (}\FloatTok{0.0}\OperatorTok{,} \FloatTok{100.0}\NormalTok{)}\OperatorTok{,}\NormalTok{ (}\FloatTok{0.0}\OperatorTok{,} \FloatTok{120.0}\NormalTok{))}

\CommentTok{\# Classify under high LIF conditions}
\NormalTok{println(}\StringTok{"=== High LIF (L=150) — Pluripotent conditions ==="}\NormalTok{)}
\NormalTok{result\_high }\OperatorTok{=}\NormalTok{ classify\_system(ds\_high\_LIF}\OperatorTok{,}\NormalTok{ bounds}\OperatorTok{;}\NormalTok{ n\_samples}\OperatorTok{=}\FloatTok{200}\NormalTok{)}
\NormalTok{print\_classification(result\_high)}

\CommentTok{\# Classify under low LIF conditions  }
\NormalTok{println(}\StringTok{"}\SpecialCharTok{\textbackslash{}n}\StringTok{=== Low LIF (L=10) — Differentiation conditions ==="}\NormalTok{)}
\NormalTok{result\_low }\OperatorTok{=}\NormalTok{ classify\_system(ds\_low\_LIF}\OperatorTok{,}\NormalTok{ bounds}\OperatorTok{;}\NormalTok{ n\_samples}\OperatorTok{=}\FloatTok{200}\NormalTok{)}
\NormalTok{print\_classification(result\_low)}

\CommentTok{\# Analyse fixed points (cell states)}
\NormalTok{println(}\StringTok{"}\SpecialCharTok{\textbackslash{}n}\StringTok{=== Fixed Point Analysis ==="}\NormalTok{)}
\KeywordTok{for}\NormalTok{ (i}\OperatorTok{,}\NormalTok{ fp) }\KeywordTok{in}\NormalTok{ enumerate(result\_high.fixed\_points)}
\NormalTok{    N}\OperatorTok{,}\NormalTok{ O}\OperatorTok{,}\NormalTok{ F}\OperatorTok{,}\NormalTok{ G }\OperatorTok{=}\NormalTok{ fp.location}
    \KeywordTok{if}\NormalTok{ N }\OperatorTok{\textgreater{}} \FloatTok{50} \OperatorTok{\&\&}\NormalTok{ G }\OperatorTok{\textless{}} \FloatTok{20}
\NormalTok{        state }\OperatorTok{=} \StringTok{"Pluripotent (stem cell)"}
    \KeywordTok{elseif}\NormalTok{ G }\OperatorTok{\textgreater{}} \FloatTok{50} \OperatorTok{\&\&}\NormalTok{ N }\OperatorTok{\textless{}} \FloatTok{20}
\NormalTok{        state }\OperatorTok{=} \StringTok{"Differentiated"}
    \KeywordTok{else}
\NormalTok{        state }\OperatorTok{=} \StringTok{"Transition state"}
    \KeywordTok{end}
\NormalTok{    println(}\StringTok{"State $i: $state"}\NormalTok{)}
\NormalTok{    println(}\StringTok{"  Location: N=$(round(N, digits=1)), O=$(round(O, digits=1)), "} \OperatorTok{*}
            \StringTok{"F=$(round(F, digits=1)), G=$(round(G, digits=1))"}\NormalTok{)}
\NormalTok{    println(}\StringTok{"  Type: $(fp.type)"}\NormalTok{)}
\KeywordTok{end}

\CommentTok{\# Check for non{-}gradient (curl) dynamics}
\CommentTok{\# The paper notes that curl dynamics are ubiquitous in gene regulatory networks}
\NormalTok{println(}\StringTok{"}\SpecialCharTok{\textbackslash{}n}\StringTok{=== Curl Analysis ==="}\NormalTok{)}
\NormalTok{test\_point }\OperatorTok{=}\NormalTok{ [}\FloatTok{60.0}\OperatorTok{,} \FloatTok{50.0}\OperatorTok{,} \FloatTok{40.0}\OperatorTok{,} \FloatTok{20.0}\NormalTok{]  }\CommentTok{\# Near pluripotent state}
\NormalTok{curl }\OperatorTok{=}\NormalTok{ curl\_magnitude(ds\_high\_LIF}\OperatorTok{,}\NormalTok{ test\_point)}
\NormalTok{ratio }\OperatorTok{=}\NormalTok{ curl\_to\_gradient\_ratio(ds\_high\_LIF}\OperatorTok{,}\NormalTok{ test\_point)}
\NormalTok{println(}\StringTok{"Curl magnitude at test point: $(round(curl, digits=4))"}\NormalTok{)}
\NormalTok{println(}\StringTok{"Curl/gradient ratio: $(round(ratio, digits=4))"}\NormalTok{)}

\KeywordTok{if}\NormalTok{ ratio }\OperatorTok{\textgreater{}} \FloatTok{0.1}
\NormalTok{    println(}\StringTok{"→ Significant non{-}gradient dynamics present"}\NormalTok{)}
\NormalTok{    println(}\StringTok{"  Forward and reverse differentiation paths will differ (see paper Fig 6)"}\NormalTok{)}
\KeywordTok{end}
\end{Highlighting}
\end{Shaded}

\textbf{Expected output:}

The stem cell model exhibits: 1. \textbf{Multiple stable states}:
Pluripotent (high Nanog, low Gata6) and differentiated (low Nanog, high
Gata6) 2. \textbf{Non-zero curl}: The system is not a pure gradient
system, meaning minimum action paths differ for differentiation vs
reprogramming 3. \textbf{Transition state}: An unstable fixed point
between the two stable states 4. \textbf{LIF-dependent landscape}:
Changing LIF concentration reshapes the potential landscape

This connects to the paper's key insight: the presence of curl dynamics
means that observing differentiation trajectories does not directly
reveal reprogramming paths.

\hypertarget{example-5-analysing-landscape-structure}{%
\subsubsection{Example 5: Analysing Landscape
Structure}\label{example-5-analysing-landscape-structure}}

\begin{Shaded}
\begin{Highlighting}[]
\KeywordTok{using}\NormalTok{ FlowClass}

\CommentTok{\# Simple bistable system (toggle switch)}
\KeywordTok{function}\NormalTok{ toggle\_switch(x}\OperatorTok{;}\NormalTok{ a}\OperatorTok{=}\FloatTok{1.0}\OperatorTok{,}\NormalTok{ n}\OperatorTok{=}\FloatTok{2}\NormalTok{)}
    \KeywordTok{return}\NormalTok{ [}
\NormalTok{        a }\OperatorTok{/}\NormalTok{ (}\FloatTok{1} \OperatorTok{+}\NormalTok{ x[}\FloatTok{2}\NormalTok{]}\OperatorTok{\^{}}\NormalTok{n) }\OperatorTok{{-}}\NormalTok{ x[}\FloatTok{1}\NormalTok{]}\OperatorTok{,}
\NormalTok{        a }\OperatorTok{/}\NormalTok{ (}\FloatTok{1} \OperatorTok{+}\NormalTok{ x[}\FloatTok{1}\NormalTok{]}\OperatorTok{\^{}}\NormalTok{n) }\OperatorTok{{-}}\NormalTok{ x[}\FloatTok{2}\NormalTok{]}
\NormalTok{    ]}
\KeywordTok{end}

\NormalTok{ds }\OperatorTok{=}\NormalTok{ DynamicalSystem(toggle\_switch}\OperatorTok{,} \FloatTok{2}\NormalTok{)}
\NormalTok{bounds }\OperatorTok{=}\NormalTok{ ((}\FloatTok{0.0}\OperatorTok{,} \FloatTok{2.0}\NormalTok{)}\OperatorTok{,}\NormalTok{ (}\FloatTok{0.0}\OperatorTok{,} \FloatTok{2.0}\NormalTok{))}

\NormalTok{result }\OperatorTok{=}\NormalTok{ classify\_system(ds}\OperatorTok{,}\NormalTok{ bounds)}

\CommentTok{\# Examine fixed points}
\KeywordTok{for}\NormalTok{ fp }\KeywordTok{in}\NormalTok{ result.fixed\_points}
\NormalTok{    println(}\StringTok{"Fixed point at $(round.(fp.location, digits=3))"}\NormalTok{)}
\NormalTok{    println(}\StringTok{"  Type: $(fp.type)"}\NormalTok{)}
\NormalTok{    println(}\StringTok{"  Eigenvalues: $(round.(fp.eigenvalues, digits=3))"}\NormalTok{)}
    
    \KeywordTok{if}\NormalTok{ fp.}\KeywordTok{type} \OperatorTok{==}\NormalTok{ SADDLE}
\NormalTok{        println(}\StringTok{"  → This is a transition state between cell fates"}\NormalTok{)}
    \KeywordTok{end}
\KeywordTok{end}

\CommentTok{\# Check landscape representation}
\NormalTok{can\_rep}\OperatorTok{,} \KeywordTok{type}\OperatorTok{,}\NormalTok{ desc }\OperatorTok{=}\NormalTok{ has\_landscape\_representation(result)}
\NormalTok{println(}\StringTok{"}\SpecialCharTok{\textbackslash{}n}\StringTok{Landscape: $desc"}\NormalTok{)}
\end{Highlighting}
\end{Shaded}

\hypertarget{theoretical-background}{%
\subsection{Theoretical Background}\label{theoretical-background}}

\hypertarget{gradient-systems}{%
\subsubsection{Gradient Systems}\label{gradient-systems}}

A gradient system satisfies dx/dt = −∇V(x) for some scalar potential
V(x). Key properties:

\begin{itemize}
\tightlist
\item
  \textbf{Jacobian symmetry}: J = −H(V) where H is the Hessian, so J =
  Jᵀ
\item
  \textbf{Curl-free}: ∇ × \textbf{F} = 0 (in 3D) or more generally, the
  Jacobian is symmetric
\item
  \textbf{No periodic orbits}: Trajectories always descend the potential
\item
  \textbf{Path independence}: Line integrals are path-independent
\end{itemize}

The condition ∂fᵢ/∂xⱼ = ∂fⱼ/∂xᵢ is both necessary and sufficient for the
existence of a potential.

\hypertarget{gradient-like-systems}{%
\subsubsection{Gradient-like Systems}\label{gradient-like-systems}}

Gradient-like systems possess a global Lyapunov function but may have
non-symmetric Jacobians away from fixed points. They share the key
property that trajectories cannot form closed loops.

\hypertarget{morse-smale-systems}{%
\subsubsection{Morse-Smale Systems}\label{morse-smale-systems}}

Morse-Smale systems allow hyperbolic periodic orbits (limit cycles)
while maintaining structural stability. They require:

\begin{enumerate}
\def\labelenumi{\arabic{enumi}.}
\tightlist
\item
  Finitely many hyperbolic fixed points
\item
  Finitely many hyperbolic periodic orbits
\item
  Transverse intersection of stable/unstable manifolds
\item
  No non-wandering points other than fixed points and periodic orbits
\end{enumerate}

\hypertarget{non-gradient-dynamics-and-curl}{%
\subsubsection{Non-Gradient Dynamics and
Curl}\label{non-gradient-dynamics-and-curl}}

As discussed by Brackston et al.~(2018), most biological systems exhibit
non-gradient dynamics. The vector field can be decomposed as:

\[\mathbf{F}(\mathbf{x}) = -\nabla U(\mathbf{x}) + \mathbf{F}_U(\mathbf{x})\]

where U is the potential (related to the probability landscape) and
\textbf{F}\_U is the curl/flux component. The curl component:

\begin{itemize}
\tightlist
\item
  Is indicative of non-equilibrium dynamics
\item
  Causes forward and reverse transition paths to differ
\item
  Cannot be inferred from static snapshot data alone
\item
  Arises naturally in gene regulatory networks due to asymmetric
  interactions
\end{itemize}

\hypertarget{connection-to-waddingtons-landscape}{%
\subsubsection{Connection to Waddington's
Landscape}\label{connection-to-waddingtons-landscape}}

The classification hierarchy relates directly to interpretations of
Waddington's epigenetic landscape:

\begin{longtable}[]{@{}ll@{}}
\toprule
System Class & Landscape Interpretation \\
\midrule
\endhead
Gradient & True potential landscape; elevation = −log(probability) \\
Gradient-like & Quasi-potential exists; landscape approximation valid \\
Morse-Smale & Local potentials around attractors; limit cycles as
valleys \\
General & Landscape metaphor breaks down; curl dynamics dominate \\
\bottomrule
\end{longtable}

\hypertarget{dependencies}{%
\subsection{Dependencies}\label{dependencies}}

\begin{itemize}
\tightlist
\item
  \href{https://github.com/JuliaDiff/ForwardDiff.jl}{ForwardDiff.jl} ---
  Automatic differentiation for Jacobians
\item
  \href{https://github.com/JuliaNLSolvers/NLsolve.jl}{NLsolve.jl} ---
  Nonlinear equation solving for fixed points
\item
  \href{https://github.com/SciML/OrdinaryDiffEq.jl}{OrdinaryDiffEq.jl}
  --- ODE integration for trajectories and manifolds
\item
  \href{https://docs.julialang.org/en/v1/stdlib/LinearAlgebra/}{LinearAlgebra}
  --- Standard library
\end{itemize}

\hypertarget{contributing}{%
\subsection{Contributing}\label{contributing}}

Contributions are welcome! Please feel free to submit issues and pull
requests.

%

\hypertarget{licence}{%
\subsection{Licence}\label{licence}}

MIT License --- see \url{LICENSE} for details.

\end{document}